# Primes, $\pi$, and Irrationality Measure

## Jonathan Sondow

ABSTRACT. A folklore proof of Euclid's theorem on the infinitude of primes uses the Euler product and the irrationality of $\zeta(2) = \pi^2/6$. A quantified form of Euclid's Theorem is Bertrand's postulate $p_{n+1} < 2p_n$. By quantifying the folklore proof using an irrationality measure for $6/\pi^2$, we give a proof (communicated to Paulo Ribenboim in 2005) of a much weaker upper bound on $p_{n+1}$.

**1. INTRODUCTION.** A folklore proof of Euclid's theorem on the infinitude of primes uses the Euler product for the zeta function [**HW**, p. 246]

$$\zeta(s) := \sum_{n=1}^{\infty} \frac{1}{n^s} = \prod_{p \text{ prime}} \frac{1}{1-p^{-s}} \qquad (\Re(s) > 1)$$

together with Euler's formula $\zeta(2) = \pi^2/6$ [**BB**, p. 383] and the irrationality of $\pi^2$ [**BB**, p. 353], [**HW**, p. 47]. If there were only finitely many primes, then $\zeta(2)$ would be rational, a contradiction.

A quantified version of Euclid's theorem is Bertrand's postulate [**HW**, p. 343]: $p_{n+1} < 2p_n$ for all $n \geq 1$, where $p_n$ denotes the $n$th prime.

In this note we quantify the folklore proof, using the notion of irrationality measure [**BB**, Section 11.3].

**Definition 1**. Let $\xi$ be an irrational number. A positive real number $\mu$ is an *irrationality measure* for $\xi$ if

$$\left| \xi - \frac{a}{b} \right| > \frac{1}{b^\mu}$$

for all integers $a$ and $b$ with $b$ sufficiently large. (Note that $\mu > 2$, because $\xi$ always has infinitely many rational approximations $a/b$ with $\left| \xi - \frac{a}{b} \right| < b^{-2}$ [**HW**, Sections 11.1, 11.3].)

We prove the following result, which gives a much weaker upper bound on $p_{n+1}$ than Bertrand's postulate.

**Theorem 1**. *If $\mu$ is an irrationality measure for $6/\pi^2$, then for all sufficiently large $n$*

$$p_{n+1} < (p_1 p_2 \cdots p_n)^{2\mu}.$$

We communicated Theorem 1 and its proof to Paulo Ribenboim in 2005 [**S**]. Recently, irrationality measures for $\zeta(2)$ have been used to prove other results on the distribution of primes [**K**], [**MSW**].

**2. PROOF OF THEOREM 1.** For $n = 1, 2, \ldots$, let

$$\frac{a_n}{b_n} := \prod_{k=1}^{n}(1 - p_k^{-2}) = \prod_{k=1}^{n}\frac{p_k^2 - 1}{p_k^2}$$

in lowest terms. Using the Euler product, we see that $6/\pi^2 = \zeta(2)^{-1} < a_n/b_n < 1$. Then, by the hypothesis and Definition 1, for $n$ sufficiently large we have

$$\frac{1}{b_n^\mu} < \frac{a_n}{b_n} - \frac{6}{\pi^2} = \frac{a_n}{b_n}\left(1 - \prod_{k=n+1}^{\infty}\frac{p_k^2 - 1}{p_k^2}\right)$$

$$< 1 - \prod_{k=n+1}^{\infty}\left(1 - \frac{1}{p_k^2}\right) < \sum_{k=n+1}^{\infty}\frac{1}{p_k^2} < \frac{1}{p_{n+1}}.$$

Since $b_n \leq (p_1 p_2 \cdots p_n)^2$, the theorem follows. ●

*209 West 97th Street, New York, NY 10025          jsondow@alumni.princeton.edu*